\documentclass[11pt]{amsart}

\usepackage{amsmath,amssymb}
\usepackage{amsthm}
\usepackage[nobysame]{amsrefs}
\usepackage{indentfirst}
\usepackage{mathrsfs}
\usepackage{xcolor}
\usepackage{amsfonts}
\usepackage{graphicx}
\usepackage[margin=1in]{geometry}
\usepackage{caption,subcaption}
\usepackage[norelsize,linesnumbered,vlined,ruled,algo2e]{algorithm2e}
\usepackage{enumerate}

\usepackage{hyperref}

\graphicspath{ {./fig1D/} }

\theoremstyle{plain}
\newtheorem{theorem}{Theorem}

\theoremstyle{definition}

\newtheorem{algorithm}[theorem]{Algorithm}

\theoremstyle{remark}

\DeclareMathOperator{\Var}{Var}

 \newcommand{\figref}[1]{Figure~\ref{#1}}

\newcommand{\ie}{\textit{i.e.}}

\newcommand{\ud}{\,\mathrm{d}}
\newcommand{\rd}{\mathrm{d}}

\newcommand{\RR}{\mathbb{R}}

\newcommand{\Or}{\mathcal{O}}

\newcommand{\EE}{\mathbb{E}}
\newcommand{\E}{\mathbb{E}}

\newcommand{\wh}[1]{\widehat{#1}}

\newcommand{\abs}[1]{\lvert#1\rvert}

\newcommand{\CC}{\mathcal{C}}
\newcommand{\cl}{R}
\newcommand{\dom}{D}

\newcommand{\Real}{\mathbb{R}}

\newcommand{\barint}{\kern4pt \raise3.4pt\hbox{\vrule height.6pt
    width7pt} \kern-11pt \int}

\begin{document}

\title{Efficient rare event simulation for failure problems in random
  media}

\author{Jingchen Liu}
\address{Department of Statistics, Columbia University, New York, NY}
\email{jcliu@stat.columbia.edu}

\author{Jianfeng Lu}

\address{Department of Mathematics, Physics, and Chemistry, Duke
  University, Box 90320, Durham, NC 27708 USA}
\email{jianfeng@math.duke.edu}

\author{Xiang Zhou}
\address{Department of Mathematics, City University of Hong Kong, Hong
  Kong SAR}

\email{xizhou@cityu.edu.hk}

\date{\today}

\begin{abstract}
  In this paper we study rare events associated to solutions of
  elliptic partial differential equations with spatially varying
  random coefficients.  The random coefficients follow the lognormal
  distribution, which is determined by a Gaussian process. This model
  is employed to study the failure problem of elastic materials in
  random media in which the failure is characterized by that the
  strain field exceeds a high threshold.  We propose an efficient
  importance sampling scheme to compute small failure probabilities in
  the high threshold limit.  The change of measure in our scheme is
  parametrized by two density functions.  The efficiency of the
  importance sampling scheme is validated by numerical examples.
\end{abstract}

\maketitle


\section{Introduction}

The study and computation of rare events in stochastic systems have
received intensive attention in recent years. Rare events, though do
not occur often, represent the most severe consequence of uncertainty
and random effects. The study of these rare events hence gives crucial
understanding and has important applications. However, due
to the small probability of occurrence of such events, the
quantification casts a serious challenge for conventional probabilistic
methods. For example, a direct Monte Carlo strategy to estimate the vanishing small probability will require a huge number of sample points to give estimates with small relative error; in other words, the huge relative variance of these estimators make them incapable of accurate prediction.

In this work, we aim at developing an efficient important sampling
strategy to study rare events associated with materials failure
problem. The method we develop in this work applies to the general
linear elasticity model for the materials failure problem.  For
simplicity, we will restrict our discussions here to a scalar model in
two dimension, which can be viewed as a model for out-of-plane
deformation of an elastic membrane under external forcing. Similar
equations also arise from other contexts, such as groundwater
hydraulics, electrostatic response of a planar media, etc.  Let $\dom
\subset \RR^2$ be an open domain with smooth boundary, which is the
equilibrium configuration of the membrane. We consider out-of-plane
displacement field $u$ given by the following boundary value problem
\begin{equation}\label{eq:PDE}
  \begin{cases}
    - \nabla \cdot (a(x) \nabla u (x)) = f(x) & \text{for  } x\in\dom;  \\
    u (x)= 0 & \text{for } x\in\partial \dom.
  \end{cases}
\end{equation}
Here $f$ is the body force acting on the material and $a: D \to \RR$
gives the stiffness of the material. We assume that the membrane is
attached to a frame at the boundary $\partial D$ and hence the
Dirichlet boundary condition in \eqref{eq:PDE}. We assume that the
external force $f$ is bounded, that is, there exist a constant $C\in
\Real$ such that
\begin{equation}
  \label{eq:fbound}
  \abs{f(x)} \leq C , \qquad x\in\dom.
\end{equation}

We study the behavior of the material under the influence of internal
randomness, which may be a result of material processing or the
uncertainty of the material properties at the microscopic level. We
adopt a probabilistic viewpoint of the complexity and heterogeneity
inherent in the material; and hence view the stiffness coefficient
$a(x)$ as a positive random field. To be more specific, we assume that
$a(x)$ is a lognormal random field, that is,
\begin{equation}\label{a}
  a(x) = \exp(- \xi(x)), \qquad x \in \dom,
\end{equation}
where $\xi$ is a stationary Gaussian random field. 
The lognormal assumption is often used in failure modeling, as it yields good
fittings to data (see for example \cite{Ohring:98}). It is also quite
natural from the mathematical point of view, as the equation is then
almost surely uniformly elliptic.
To simplify notation and without loss of generality, we assume that $\E \xi(x) = 0 $ and $\Var \xi(x) =1$.

The random field viewpoint is taken in the homogenization theory for
random heterogeneous and composite materials (see
e.g. \cites{Milton:02, Torquato:01}). However, for the study of rare
material failure events, standard homogenization theory is not enough
to capture rare events, despite the recent advances in the
understanding of variance scaling and central limit theorem
\cites{Bal:08, GloriaOtto:11, GloriaOtto:13, Nolen:13+}. Here our
focus is on developing efficient numerical methods for the computation
of the material failure probability via importance sampling.

We consider materials failure such that under the external force, the
strain of the elastic deformation exceeds a prefixed level at some
point. More precisely, let $b \gg 1$ be the given threshold, the
failure probability is defined as
\begin{equation}\label{failprob}
  P\Bigl\{  \sup_{x\in\dom} \abs{\nabla u} \geq b  \Bigr\}.
\end{equation}
Note that this probability is very small when $b$ is large. Direct
Monte Carlo simulation is hence difficult to compute such
probabilities and achieve small relative error. Instead, we employ
importance sampling techniques for the computation of the failure
probability.

In this paper, we propose an efficient Monte Carlo method via importance sampling to compute
the small failure probabilities as in \eqref{failprob} when the
differential equation is driven by a Gaussian random field as in
\eqref{a}.  The change of measure proposed in this paper is not of the
exponential tilting form and therefore is nonstandard.  In the
one-dimensional setting, the algorithm can be proved to be
asymptotically efficient. For the case in higher dimensions, due to
the lack of large deviations results, efficiency cannot be rigorously
established. However, the algorithm does admit a very good performance in our numerical studies.

\medskip

It is a long history of studies of material failure and structure
safety in the civil engineering and material sciences from the
probabilistic viewpoint or the extreme value theory.  The important
role of Weibull distribution \cite{Weibull1951} for (idealized)
weak-link principle is well established and applied in numerous
engineer applications, although it faces lots of challenges from real
material properties.  A substantial progress toward this challenge
is the well-known work of Ba\v{z}ant
 on the statistical size
effect \cite{BazantPNAS2004}.  In contrast to these engineering
statistical  models for
material/structure failure problem, we take  the 
mechanistic approach to use the 
classical linear
elasticity model, i.e., \eqref{eq:PDE}. It is admitted that  many physical processes such as the
development of fractures and cracks as well as dynamical failure scenarios  are not shown in this
model.  Yet, it is
 a good prototype model, in the balance of tractability and
complexity of modeling material failure problem.  Nevertheless, the probabilistic study 
of this model  for high excursion of strain field  is  helpful to  
shed  light on the extreme mechanical   behaviours of random elastic media.  
Furthermore, our model is quite general and is never limited to the application of elastic mechanics.
There are lots of other important physical and engineer problems 
modeled in the exactly same form of our elliptic equation \eqref{eq:PDE} with random coefficients.
For instance,  the Darcy equation with uncertain coefficients is 
 the canonical model for groundwater  study and in this context, 
the derivative of the solution  $\nabla u$  is   an important
physical quantity related to the phase speed of pollutants carried by
groundwater.


In view of  extremely small failure probabilities concerned here, 
our work   fits into the general scope of works that devote to rare
event simulations. Different approaches have
been proposed in recent year for such problems, in particular in engineering and
industrial applications. For instance, the idea of design point shift
has been used in the framework of polynomial chaos expansion for
failure events \cite{Paffrath-JCP-2007}.  The numerical adaptive
strategy is also tested on the PDE with random input data
\cite{Chen2013233}.  In terms of Monte Carlo importance sampling
method which is free of ``curse of dimensionality'', the work
\cite{LiLiXiu:11} combines the cross-entropy method and the surrogate
model to efficiently calculate the failure probabilities, which in
principle works for very general problems. 
The recent work \cite{Papanicolaou2012} applied the large deviation and importance
sampling method to the calculation of the failure probability of hypersonic engines.  We
also mention the study of rare events in optical pulses modeled by
randomly perturbed one dimensional non-linear Schr\"{o}dinger equation
(see e.g.~\cites{Donovan2011,MooreKathSIAMReview}). We refer to
\cites{DUPWAN04,BUC90,ASMGLY07,EVa10} for general techniques of
importance sampling and rare events simulations.

Our problem is also closely connected to the probabilistic theories
for Gaussian random field.  For stochastic systems driven by
light-tailed random variables (such as, Gaussian random variables), it
is customary to consider exponential change of measure for the design
of efficient importance sampling algorithm.  The parameters are
usually selected by the minimal cross-entropy method \cites{RuKr04} or
some control problems related to large deviation principle
\cite{DUPWAN04}.  For heavy-tailed stochastic systems, some recent
works are \cites{BGL07,BL08,BL12}.  In the context of Gaussian
processes and random fields, the most well studied events are the high
level excursions (tail events of the supremum) \cite{ABL09}; the tail
events of other convex functionals of Gaussian random fields are also
of interest \cites{LiuXudensity,Liu10,LiuXu11,LiuXu14}.
The method in this paper is in part built on the results in this literature.

The rest of the paper is organized as follows. The description of the
algorithm is given in Section~\ref{sec:alg}. Implementation details
and numerical results are discussed in
Section~\ref{sec:numer}. Conclusions and discussions are summarized in
Section \ref{sec:conclude}.

\section{The main method}\label{sec:alg}

\subsection{Rare event simulation, variance reduction, and importance
  sampling}

Let us consider the problem of estimating a small probability $w=P(B)
\ll 1$ or a family of probabilities $w(b) = P(B_b)$, where $b$ is the
rarity parameter that indicates the difficulty of the problem. In our
case, we identify the rarity parameter as the threshold, and the event
is given by $B_b = \{ \sup_{\dom} \abs{\nabla u} \geq b\}$.  As $b$
tends to infinity, the probability of interest $w(b) = P(B_b)$ tends
to $0$.

As the probabilities to be estimated are very small, the computational
error needs to be quantified relative to the probabilities of interest.
Let us consider an unbiased Monte Carlo estimator $Z_b$ for $w(b)$
such that $\EE Z_b = w(b)$. The relative error is given by 
$\Var(Z_b) / w^2(b)$ or $\EE(Z_b^2)/w^2(b) = 1+\Var(Z_b) / w^2(b)$.
Suppose that  $n$ independent and identically distributed replicates of $Z_b$ are generated, denoted by $Z_b^{(1)}, \ldots, Z_b^{(n)}$. Let
\begin{equation*}
\bar Z_n \triangleq \frac 1 n \sum_{i=1}^n Z_b^{(i)}
\end{equation*}
be the averaged estimator, whose variance is
\begin{equation*}
\Var (\bar Z_{n} )= \frac{\Var (Z_b)}{n}.
\end{equation*}
Via Chebyshev's inequality, we obtain that
\begin{equation*}
P(|\bar Z_n -w(b)| >\varepsilon w(b)) \leq \frac{\Var (Z_b)}{n \varepsilon^2 w^2(b)}.
\end{equation*}
For any $\delta>0$ if we intend to estimate $w(b)$ with at most $\varepsilon$ relative error with at least probability $1-\delta$, then it suffices to generate
\begin{equation*}
  n = \frac{\Var(Z_b)}{ w^2(b)} \delta^{-1} \varepsilon^{-2}
\end{equation*}
samples. Hence, the necessary sample size is proportional to the relative error ${\Var(Z_b)}/{ w^2(b)}$.

Consider the direct Monte Carlo estimator $I_{B_b}$. Its relative error is
$$\frac{\Var(I_{B_b})} {w^2(b)} = \frac{1-w(b)}{w(b)} \to \infty$$
as $w(b) \to 0$. The necessary sample size is $n = w^{-1}(b) \delta^{-1}\varepsilon^{-2}$ and $I_{B_b}$ is considered an inefficient estimator for $w(b)$ when it is very small.
There are several efficiency criteria in the literature \cites{BUC90, ASMGLY07}. 
The most widely used is the \emph{weak efficiency} or \emph{asymptotic efficiency} requiring that for all $\epsilon >0$,  $\EE(I^2_{B_b})/w^2(b) = o(w^{-\epsilon}(b))$ as $w(b)\to 0$.
In the numerical analysis, we will investigate the relative errors of the proposed estimator. The empirical study shows that our estimator admits reasonably small relative errors (less then 10) when $w(b)$ is very small (less than $10^{-6}$), although rigorous efficiency is difficult to establish and is beyond the scope of this paper.

In the subsequent analysis, we employ importance sampling as the main variance reduction technique
that is based on the following identity
\begin{equation*}
P(B) = \int I(\omega \in B) P(\rd\omega) = \int I(\omega \in B) \frac{\rd P}{\rd Q} Q(\rd\omega), \quad \mbox{for all measurable set $B$.}
\end{equation*}
for a measure $Q$ such that $Q(\cdot\cap B)$ is absolutely continuous
with respect to the measure $P(\cdot \cap B)$. As a consequence, we have
\begin{equation*}
Z(\omega) = I(\omega \in B) \frac{\rd P}{\rd Q}(\omega)
\end{equation*}
is an unbiased estimator of $P(B)$ under $Q$, in other words, $\EE^QZ
= P(B)$, where we use $\EE^Q$ for the expectation with respect to the
measure $Q$,

It is easy to verify that if we choose $\mathcal Q = P(\cdot | B) =
P(\cdot \cap B)/ P(B)$ then the corresponding importance sampling
estimator admits zero variance. Thus, we often call $\mathcal Q$ the
zero-variance change of measure. On the other hand, $\mathcal Q$ is
clearly of no practical value, in that the likelihood ratio is almost
surely $P(B)$ that is precisely the quantity we want to
compute. Nevertheless, the measure $\mathcal Q$ provides a guideline
to construct a change of measure for the efficient computation of
$P(B)$. We need to construct a measure $Q$ that is close to $\mathcal
Q$ such that we are capable of sampling from $Q$ and computing the
Radon-Nikodym $\rd Q/\rd P$.

\subsection{The change of measure}

Let us characterize a measure $Q$ for the random field $\xi$ as in
\eqref{a} defined on the continuous sample path space $\CC(\dom)$,
where $\xi: \dom \to \RR$ is a realization of the random field and
$\dom$ is the domain of the PDE. Our choice of $Q$ depends on two
probability density functions $h(\cdot)$ and $g_x(\cdot)$ to be
determined later.  Given $h$ and $g$, the Radon-Nikodym derivative of
$Q$ with respect to $P$ is given by
\begin{equation}\label{LR}
  \frac{\rd Q}{\rd P}(\xi)
  = \int_{\dom} h(x) \frac{g_x(\xi(x))}{\phi_x(\xi(x))} \ud x
\end{equation}
 Here, for each point $x \in \dom$, $\phi_x(\cdot)$ is the marginal density function of $\xi(x)$ under $P$, that is,
\begin{equation*}
  \int_A \phi_x(y) \ud y = P(\xi(x) \in A)\quad \mbox{for any measurable set $A \subset \RR$.}
\end{equation*}
In \eqref{LR}, $h(\cdot)$ is a density function over the domain $\dom$
and $g_x(\cdot)$ for each $x \in \dom$ is a density function on
$\RR$. We will choose $h$ and $g_x$ such that the corresponding
measure $Q$ is a good approximation of $\mathcal Q$ for variance
reduction.

Let us first explain how to sample from $Q$ before discussing the
choices of $h$ and $g_x$. It consists of three steps:
\begin{algorithm}\ \label{algcont}
\begin{enumerate}
\item Sample a random index (position) $x \in \dom$ following the
  density $h(x)$;

\item Conditional on the realized $x$, sample a random number $\xi(x)$ following the density $g_x (\cdot)$;

\item Conditional on the realized $x$ and $\xi(x)$, sample
  $\xi$ on $\dom \backslash \{x\}$ from the conditional distribution
  $P\{\xi\in \cdot \mid \xi(x)\}$.
 \end{enumerate}
\end{algorithm}
It is easy to verify that the above three-step procedure is consistent
with the Radon-Nikodym derivative \eqref{LR}.  If $g_x = \phi_x$ for
all $x\in\dom$, then we have $Q=P$.  Thus, the distributions of $\xi$
under $P$ and $Q$ are different only at one random location that is
labeled by $x$ sampled from the density $h$.  This suggests if the
occurrence of the rare event $\{\sup_{x\in\dom} |\nabla u(x)| > b\}$
is mostly due to the abnormal behavior of the random field $\xi$ at
one location, our $Q$ would be a good candidate of approximating the
zero-variance change of measure $\mathcal{Q}$.  In this sense, the
distribution of the random index $x$ (\ie~$h(x)$) should be
approximately the distribution of the location where $\xi$ deviates
mostly from its original law under $P$. Furthermore, the distribution
$g_x$ characterizes how $\xi(x)$ deviates from its original law
$\phi_x$.  In what follows, we will describe in details the choices of
$g_x$ and $h$.

\subsection{The excursion level and the choice of $g_x$}
\label{ssec:exlevel}

Bearing in mind the above intuition, we proceed to describing $h$ and
$g_x$. Among these two, $g_x$ is more important as it quantifies the deviation of $\xi$ from its original law.  The basic
idea is as follows. If $\sup_{\dom} |\nabla u(x)|$ admits an
excursion over some high level $b$, then the process $\sup_{\dom}
\xi(x)$ must also have a high excursion over some level $l_x$
depending on $b$ and the precise location $x$ where the excursion
occurs.  This observation is due to the connection between $\xi$ and
$u$ in the PDE \eqref{eq:PDE}. Therefore, we expect that, under the
distribution $g_x$, $\xi(x)$ attains a high level $l_x$ that goes to
infinity at some rate with $b$.  Once the level $l_x$ has been
decided, we would choose the distribution $g_x(\cdot)$ be a Gaussian
distribution with mean $l_x$.

In the analysis of the one-dimensional equation, explicit formula is
available for the solution to the differential equation,
$$u'(x) = e^{\xi(x)} \Big\{ -F(x) + \frac{\int_\dom F(x) e^{\xi(y)}dy}{\int_\dom e^{\xi(y)}dy}\Big\}$$
where $F(x) = \int_0^x f(y) dy.$ Notice that $F(x)$ is a bounded
function and thus $\log u'(x) = \xi(x) + \Or(1)$. Hence $\max |u'(x)|
> b $ implies $\max|\xi(x)| > \log b +\Or(1)$.  Based on this closed
form solution, it is reasonable to consider that $l_x$ is
approximately $\log b$.  The optimal choice of $l_x$ would be of order
$\log b + \Or(\log\log b)$.  In the high-dimensional analysis, the PDE
does not have a closed form solution.  It is generally difficult to
derive an analytic relationship between $b$ and $l_x$. Nevertheless,
we conjecture that the relationship $l_x \approx \log b$ is generally
appropriate. This would be justified in our numerical examples.

Based on the above discussion, we choose $l_x$ such that it is just enough for $\sup_{\dom} |\nabla u(x)|$ to exceed $b$. In particular, for each $x_0\in \dom$ and $l$, we define
$$\xi_{l,x_0}(x) = l \times C(x-x_0) = E\{\xi(x) \mid  \xi(x_0) =l\}$$
where $C(\cdot)$ is the covariance function of $\xi$ that is
\begin{equation*}
  C(x) = \mbox{Cov}(\xi(y), \xi(y+x)).
\end{equation*}
The fact $l \times C(x-x_0) = E\{\xi(x) \mid  \xi(x_0) =l\}$ is due to that $\xi(x)$ has zero mean and unit variance.
For other cases, the form of conditional expectation can be adapted.
Let $u_{l,x_0}(x)$ be the solution to the PDE \eqref{eq:PDE} with $a(x) = e^{-\xi_{l,x_0}(x)}$. Then, $l_{x_0}$ is given by
\begin{equation}\label{lx}
l_{x_0}  = \min\,\{l: \max_{x\in \dom}|\nabla u_{l,x_0}(x)| \geq b\}.
\end{equation}


We now provide an intuitive explanation for the above choice of $l_x$.
The basic understanding is that the high excursion of $|\nabla u|$ is
caused by the high excursion of the input process $\xi$. In order to
determine the necessary excursion level, we perform the following
calculations.  
Conditional on $\xi(x_0) = l_{x_0}$ that is a large
number, the conditional field $\xi$ has the following representation
\begin{equation*}
  \xi(x) = l_{x_0} C(x-x_0) + r(x-x_0)
\end{equation*}
where $r(x)$ is a zero-mean Gaussian process whose covariance function
can be obtained by conditional Gaussian calculations. The rationale of
\eqref{lx} is as follows. The process $r(x)$ is the remainder process
after taking out the conditional mean and $r(x)$ is of a constant
order.  If $l_{x_0}$ is selected to be large, then the variation of
$r(x)$ is negligible compared to the conditional mean.  Therefore, the
conditional field can be approximated by
\begin{equation*}
  \xi(x)\approx l_{x_0} C(x-x_0).
\end{equation*}
By solving \eqref{lx}, $l_{x_0}$ is the minimum level  that $\xi(x)$ needs to achieve such that $\sup |\nabla u(x)|$ just exceeds $b$.

Having $l_{x}$ defined, we then choose $g_{x}$ to be the Gaussian distribution
\begin{equation}\label{gx}
g_{x} \sim \mathcal{N}(l_{x}, l_{x}^{-2}).
\end{equation}
The choice of the variance of $g_x$ comes from the following
intuition. The function $l_x$ is interpreted as the necessary level
$\xi(x)$ needs to exceed. Notice that $\xi(x)$ is marginally a
standard Gaussian random variable. Conditional on $\xi(x) > l_x$,
$\xi(x) - l_x$ is asymptotically an exponential random variable with
variance $l_x^{-2}$. Thus, the choice of variance for $g_x$ aims at
matching the scale of the overshoot of the standard Gaussian random
variable. In the simulation study, we also vary the choices of this
variance and found that $l_{x}^{-2}$ yields the best numerical results
in terms of variance reduction.  For each $x_0 \in \dom$, we provide
an iterative algorithm to compute $l_{x_0}$.

\begin{algorithm}\label{alg:findl}
Initialize $l_{x_0}^{(0)}= \log b$ and $n=0$

\While{``not converge''}{
  \begin{enumerate}
  \item Solve the PDE \eqref{eq:PDE} for $\xi(x) = l_{x_0}^{(n)} C(x -
    x_0)$ numerically. Denote the solution by $u^{(n)}$.

  \item Set $l_{x_0}^{(n+1)} = l_{x_0}^{(n)} - \log \sup\abs{\nabla
      u^{(n)}} + \log b$.
  \end{enumerate}
}
\end{algorithm}

When converged, the above algorithm yields a level $l_{x_0}^{(\infty)}$ satisfying \eqref{lx}.
	Furthermore, if $\xi(x) = l_{x_0} C(x-x_0)$, then we expect that $\sup_{\dom} |\nabla u|$ and $l_{x_0}$ have the following relationship
\begin{equation}\label{rel}
  \sup_{\dom} |\nabla u| \sim \kappa_{x_0} l_{x_0}^\alpha e^{l_{x_0}}.
\end{equation}
This relationship is correct for the one-dimensional case, and is
conjectured for the high-dimensional case (we will verify this
numerically in Section~\ref{sec:numer}). Assuming \eqref{rel},
Algorithm~\ref{alg:findl} is asymptotically the Newton-Raphson
algorithm for the equation
\begin{equation*}
  \kappa_{x_0} l_{x_0}^\alpha e^{l_{x_0}} = b.
\end{equation*}
Given the initial value $l^{(0)}_{x_0} = \log b$, it can be shown that
$|l^{(2)}_{x_0} - l^{(\infty)}_{x_0}| = o(1/\log b)$ which is accurate
enough for the construction of an efficient importance sampling.
Therefore, in our implementation, we take only very few iterations in
Algorithm~\ref{alg:findl}.

\subsection{The choice of $h$}

We now proceed to the other tuning parameter, the distribution $h$.
Recall that $x$ localizes the largest deviation from the original
distribution. Such a deviation is quantified by the level $l_x$ for
each $x$. For each $x\in \dom$, we then choose
\begin{equation*}
  h(x)  \propto P\left\{\xi(x)> l_x  \Big \vert  \sup_\dom | \nabla u(x) | > b\right\}\propto P\{\xi(x) > l_x\}.
\end{equation*}
That is, conditioning on the occurrence of the rare event, the probability that $x$ is in a small neighborhood should be proportional to the probability that $\xi(x)$ exhibits a high excursion.
After normalization, we get
\begin{equation}\label{h}
  h(x) = \frac{P(\xi(x) > l_x)}{\int_{y\in \dom} P(\xi(y) > l_y)dy} .
\end{equation}
Sampling from $h$ requires the numerical evaluation of $l_x$ which
induces some computational overhead. To further reduce the
computational complexity, we will evaluate $l_x$ for a finite grid
that spreads over the domain $\dom$ and use an interpolation for the
rest of the domain. Details will be presented in the subsequent
section.

\subsection{Summary}

Based on the above description, we generate the process $\xi(x)$ according to Algorithm \ref{algcont}. The tuning distributions $g_x$ and $h$ are given in \eqref{gx} and \eqref{h} where $l_x$ is defined through \eqref{lx}.
Finally the importance sampling estimator is 
$$Z_b= I(\sup|\nabla u(x)| > b)\Big(\int_{\dom} h(x) \frac{g_x(\xi(x))}{\phi_x(\xi(x))} \ud x\Big)^{-1}.$$
For the implementation, we discretize the space, the details of which is presented in the following section.


\section{Numerical Details and Examples}\label{sec:numer}

We will present numerical examples for one-dimensional and
two-dimension cases.  The domain for the one dimensional PDE is
$[0,1]$ and the domain for the two dimensional PDE is $\dom =
[0,1]\times [0,1]$. The covariance function of the Gaussian field
$\xi$ in both cases is
\begin{equation}\label{eqn:CovFun}
E\{\xi(x)\xi(y)\} =C(x-y)= \exp(-|x-y|^2/\cl^2),
\end{equation}
where $|\cdot|$ is the Euclidean distance.
The scalar $\cl$ is known as the correlation length of the random field.
The boundary condition is Dirichlet $u|_{\partial \dom}=0$  unless it is specified otherwise.
Our numerical results consist of a verification of the exponential relationship between $\max \xi(x)$ and $\max|\nabla u(x)|$, visualization of the excursion level $l_x$, and empirical performance of the importance sampling algorithm in estimating the failure probabilities.

\subsection{Implementation of the algorithm}

To sample the random field $\xi$, we first need to discretize the
domain $\dom $ and use the field on the discrete mesh grid as an
approximation.  Denote a point in $\RR^2$ as $x=(x^1,x^2)$.  Choose a
discrete mesh $\wh \dom_N\triangleq \{ (x^{1}_{i}, x^{2}_{j})\}$ where
$0=x^{1}_{0}< x^{1}_{1}<\ldots < x^{1}_{n_1}=1$, $0=x^{2}_{0}<
x^{2}_{1}<\ldots < x^{2}_{n_2}=1$ and $n_1\times n_2=N$.  Then any
sample path of the continuous random field $\xi(x)$ is approximated by
the $N$-dimensional random vector $\big(\xi(x^1_i,x^2_j): 1\leq i\leq
n_1, 1\leq j\leq n_2 \big)$.  The following is the discrete analogue
of Algorithm \ref{algcont} on the grid.

\begin{algorithm}\label{algdis}\
\begin{enumerate}
\item Sample a random index $\tau\in \wh \dom_N$ following the weight  $h(\tau)\delta_\tau$ as in \eqref{h}
where $\delta_\tau$ is the corresponding Lebesgue measure of the cell associated to $\tau$ in the mesh $\wh \dom_N$.
\item Conditional on the realized $\tau$, sample $\xi(\tau)$ following the density $g_\tau (\cdot)\sim N( l_\tau, l_\tau^{-2})$;
\item Conditional on the realized $\tau$ and $\xi(\tau)$, sample
  $\{\xi(x): x\in \wh \dom_N\}$ from the conditional distribution
  $P(\xi \in \cdot | \xi(\tau))$;
\item Solve PDE \eqref{eq:PDE} with $a(x)=\exp(-\xi(x))$ and calculate
  $\sup_{x\in \dom} |\nabla u(x)|$;
\item Output the estimator
  \begin{equation*}
    Z_b = I(\sup_\dom |\nabla u(x)|>b)\Big[ \sum_{x\in \wh \dom_N}   h(x)\delta_x \cdot \frac{g_x(\xi(x))}{\phi_x (\xi(x))}\Big]^{-1}.
  \end{equation*}
\end{enumerate}
\end{algorithm}

We now provide further details on steps (3) and (4) in the above
algorithm. Step (3) generates a random vector $\xi$ on $\wh \dom_N$
given a realization $\xi(\tau)$. Notice that $\xi$ on $\wh \dom_N$ is
a multivariate Gaussian random vector and thus the conditional
distribution is still multivariate Gaussian. The conditional mean is
$\EE \{\xi(x) \mid \xi(\tau)=y\} =y \times C(x-\tau)$
and the conditional covariance matrix is $\wh C_{N-1,N-1} - \wh
C_{N-1,1} \wh C_{N-1,1}^T$.  Here, $\wh C_N$ is the $N$ by $N$
unconditional covariance matrix, $\wh C_{N-1,1} $ is column of $\wh
C_N$ corresponding to $\tau$ with row corresponding to $\tau$ deleted,
and $\wh C_{N-1,N-1} $ is the $(N-1)\times(N-1)$ sub-matrix of $\wh
C_N$ with row and column corresponding to $\tau$ deleted.
However, in generating this conditional sample, we do not need
to decompose the conditional covariance matrix.
A simple procedure
${ \xi(x)} \triangleq   \xi'(x) +\wh C (x - \tau) (  \xi(\tau) -  \xi'(\tau))$ for $k=1,2,\ldots,N$, for all $x\in \wh D_N$
gives the conditional sample $\{\xi(x): x\in \wh D_N\}$ with  the known (conditional)   value of  $\xi(\tau)$,
where $\{\xi'(x): x\in \wh D_N\}$ is the $N$-dim Gaussian vector
with covariance matrix $\wh C_N$.
To sample the
multivariate Gaussian random vectors $\{\xi'(x): x\in \wh D_N\}$ , we adopt the Cholesky
decomposition of the covariance matrix that is computed via the
pivoted Cholesky factorization \cite{PCholesky} in LAPACK 3.2. Due to
stationarity, this decomposition needs to be computed only once.

Once a realization of $\xi$ is generated, the partial differential equation \eqref{eq:PDE} is solved by a standard numerical solver.
Many traditional  advanced  numerical strategies such as  adaptive mesh refinement will  improve the efficiency and the accuracy of the PDE solver.
For the numerical examples in this paper, we use the finite volume method on a uniform  mesh in $\dom$.
It is adequate to demonstrate the performance of the variance reduction with a sufficiently high resolution of the mesh grid.

In summary, our preprocessing includes a Cholesky factorization of the
conditional covariance matrix and of solving an inverse problem to
obtain the excursion level function $l$. In the sequel, we present
various numerical results related to our rare event calculation of the
failure probabilities.

\subsection{Verification of exponential relationship}
We start from a numerical verification of our conjecture about the
maximums of $\xi$ and $|\nabla u|$. Our algorithm as well as the
analysis depends on the following scaling relationship between the
input field $\xi$ and the strain field $\nabla u$, whenever either of
them has a high excursion,
\begin{equation}
\label{eqn:exp-scale}
\sup_\dom \xi      \sim \sup_\dom \log |\nabla u|.
\end{equation}

\begin{figure}[htbp]
  \begin{center}
  \begin{subfigure}[b]{0.48\textwidth}
     \includegraphics[width=\textwidth]{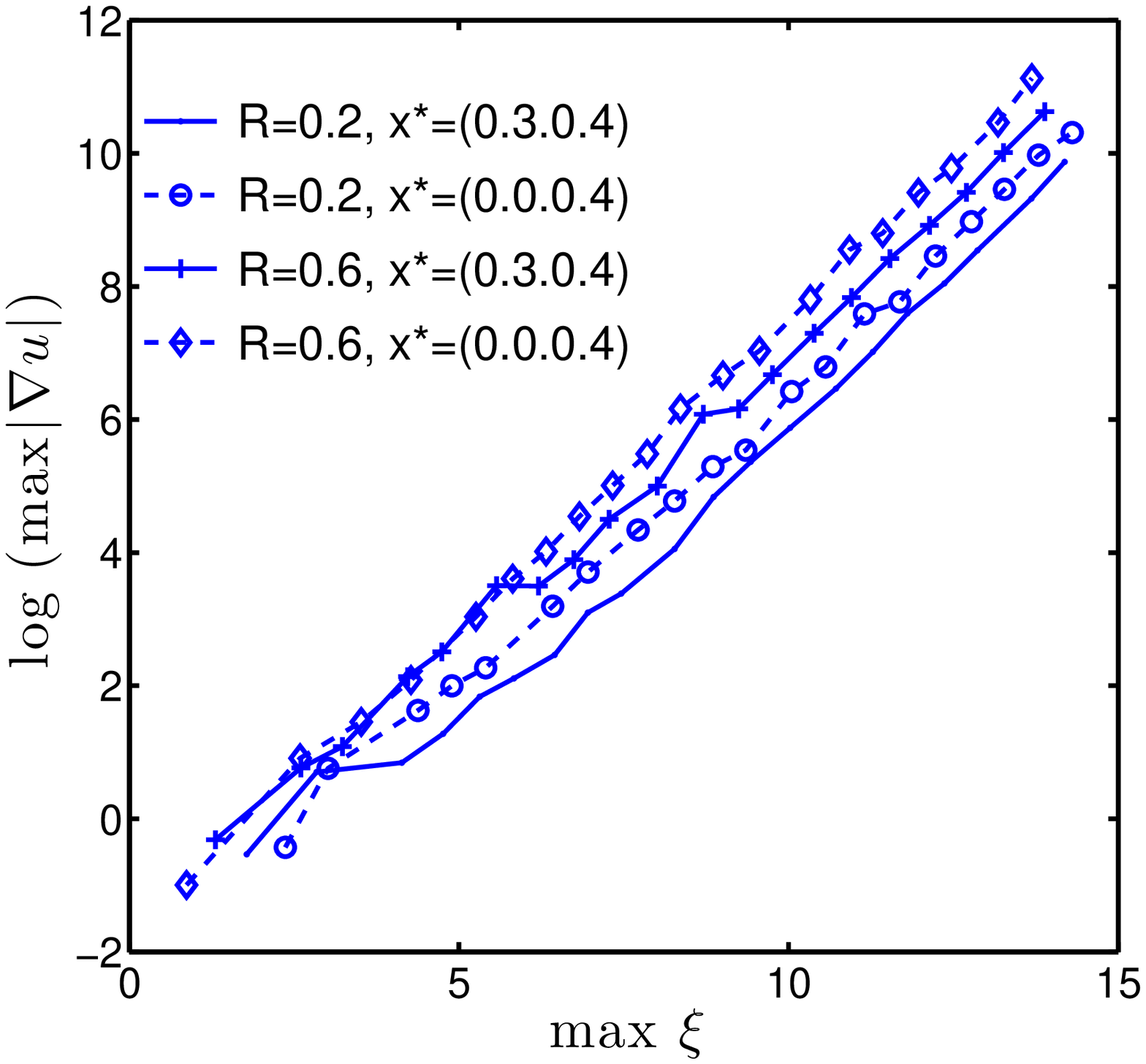}
     \caption{stochastic approach of verification}
     \label{figA:exp}
     \end{subfigure}
  \begin{subfigure}[b]{0.48\textwidth}
       \includegraphics[width=\textwidth]{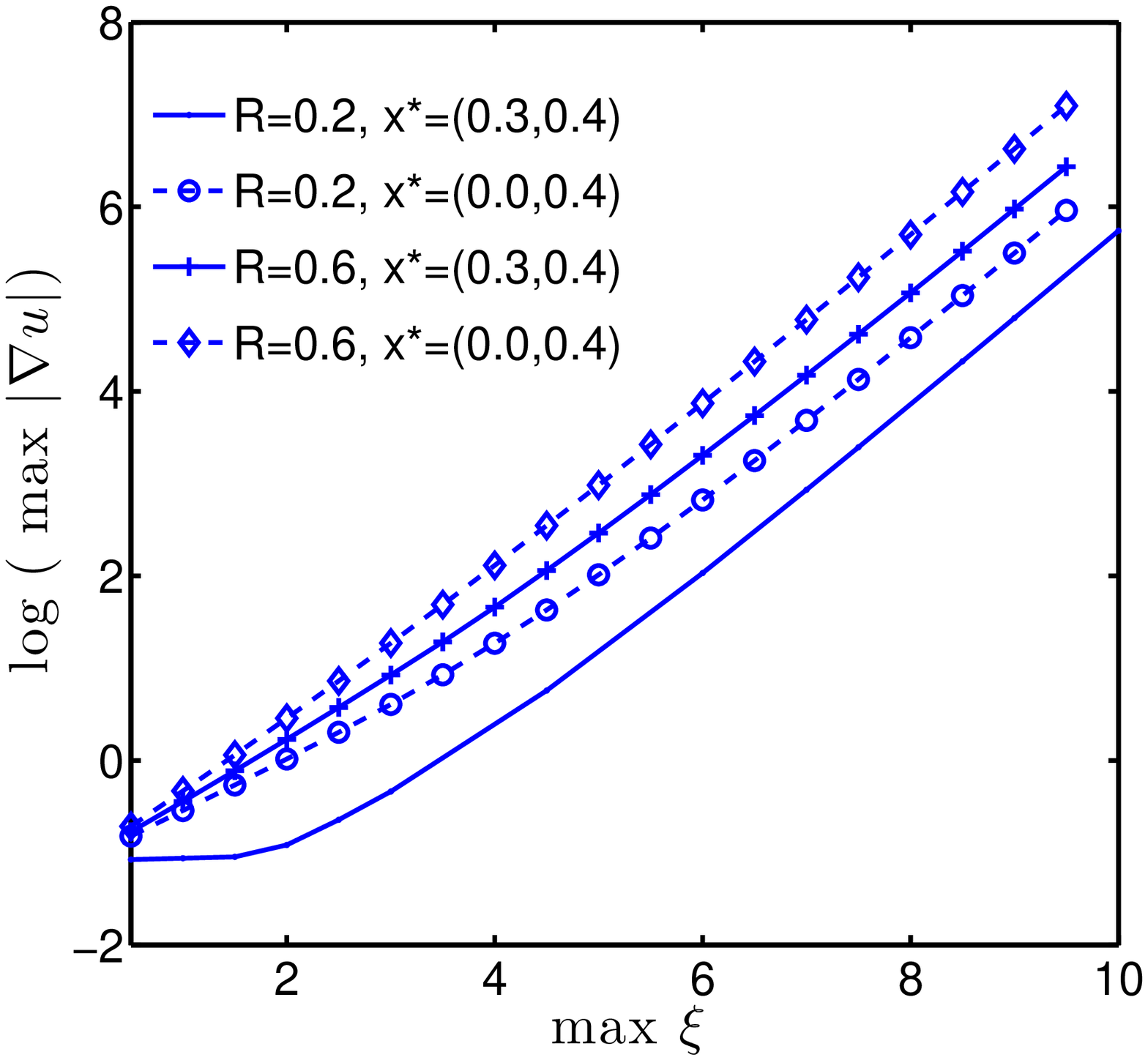}
       \caption{deterministic approach of verification}
       \label{figB:exp}
       \end{subfigure}
       \caption{Numerical verification of $\max |\nabla u| \sim
         \exp(\max(\xi))$. The domain is $\dom=[0,1]^2$.  The left
         panel is the stochastic approach where for each different
         value of $l$, one sample path of $\xi(x)$ is generated
         conditional on $\xi(x^*)=l$.  Due to the randomness, $l$ is
         not exactly equal to, but very close to, the sampled maximum
         $\max \xi$.  The right panel is the deterministic approach
         where $\xi(x)= l \times C(x-x^*)$ is deterministic and $l$ is
         precisely the maximum of $\xi$.  }
    \label{fig:exp} 
\end{center}
\end{figure}

In the one-dimensional case, the explicit formula of $u$ and $\nabla
u$ allow a verification of the above relation for any bounded external
force.  For the problem of higher dimension, as explicit solution is
not available, the rigorous justification of \eqref{eqn:exp-scale} is
more challenging.  Here, we numerically verify this relationship via a
stochastic approach and a deterministic approach.  In the stochastic
approach, one spatial location $x^*\in \dom$ is selected and a
sequence of the excursion levels $l$ are selected.  For each $l$, we
generate a random sample path $\xi(x)$ conditioned on $\xi(x^*)=l$.
We then calculate the maximum value of the strain $\max |\nabla u|$
corresponding to generated $\xi(x)$ with the homogeneous force
$f(x)\equiv 1$.  In the left panel of \figref{fig:exp}, we plot
$\log(\max |\nabla u|)$ versus $\max \xi$ for different spatial
locations $x^*$. For the deterministic approach, we observe that for
$l$ sufficiently large, the conditional field is approximately $\xi(x)
\approx l\times C(x-x^*)$, as discussed in Section~\ref{sec:alg}.
Hence, we solve the PDE with simply setting $\xi =l\times C(x-y)$
where $C$ is the covariance function \eqref{eqn:CovFun}. The numerical
results are in the right panel of \figref{fig:exp}.

The numerical results in \figref{fig:exp} confirm that $\max \log
|\nabla u|$ is asymptotically linearly proportional to $\max \xi$ and
thus justifies \eqref{eqn:exp-scale}. When the external force $f$ is
inhomogeneous, a similar relationship can be established numerically.
In addition, on comparing the results for different correlation
lengths, we found that smaller correlation length yields lower $\max
|\nabla u|$.  Thus, for smaller correlation length, to ensure that
$\max |\nabla u|$ reaches a level $b$, larger values of $\max \xi$ are
required corresponding to higher excursion level function $l_x$.

%
%

\subsection{Excursion level function $l_x$}
The excursion level function $l_x$ characterizes the spatial
distribution of the extreme values of $\xi(x)$ conditional on the
failure event.  We have explained how to find this function in Section
\ref{ssec:exlevel}.  The density function $h$ is determined by $l_x$
via \eqref{h} and the smaller value of $l_x$ implies a higher
likelihood of observing an excursion at or around $x$.

\begin{figure}[htbp]
\begin{center}
\begin{subfigure}[b]{0.48\textwidth}
\includegraphics[width=\textwidth]{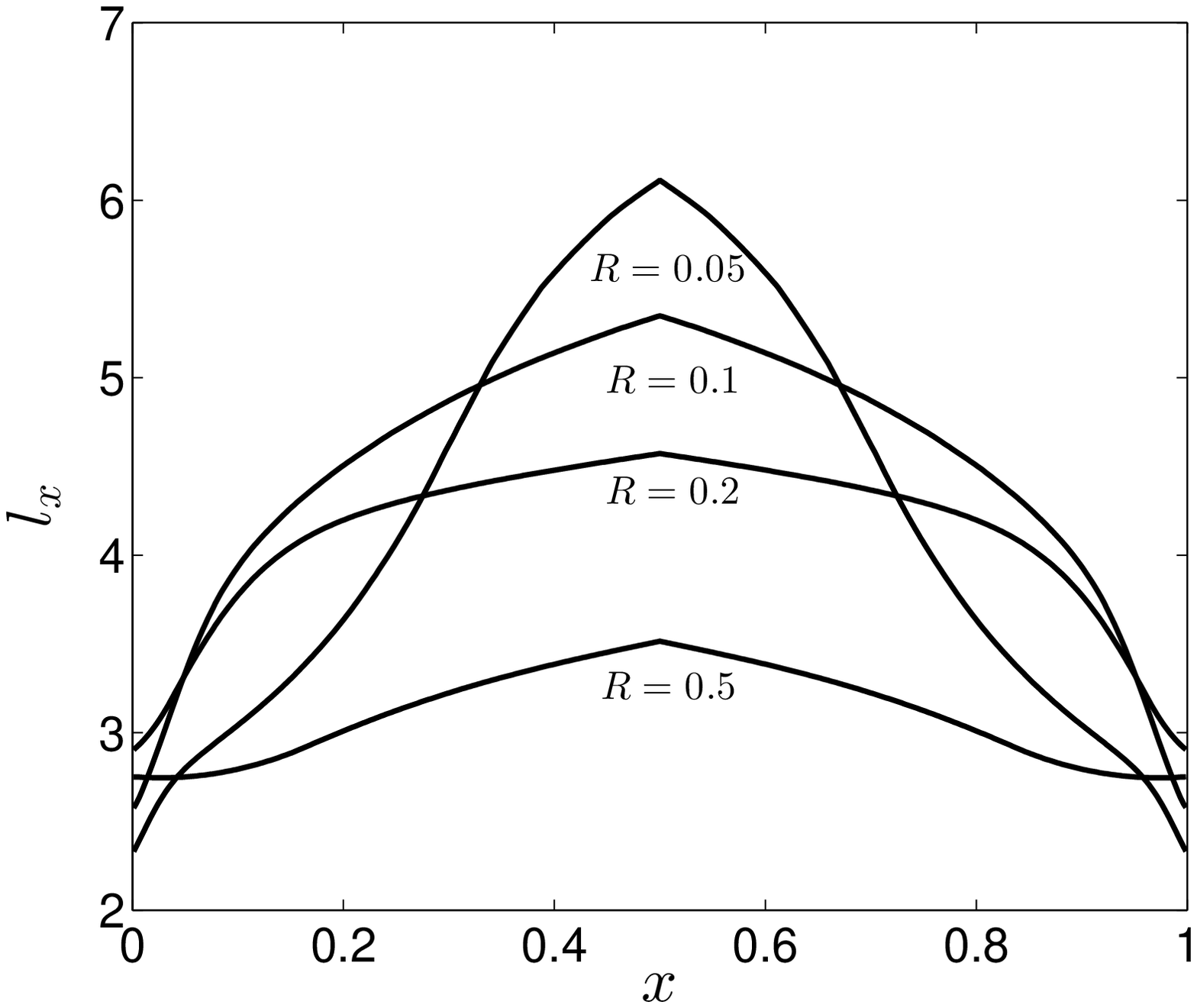}
\caption{$b=4$}
\label{fig:hfun_dc_1d-b4}
\end{subfigure}
\begin{subfigure}[b]{0.48\textwidth}
\includegraphics[width=\textwidth]{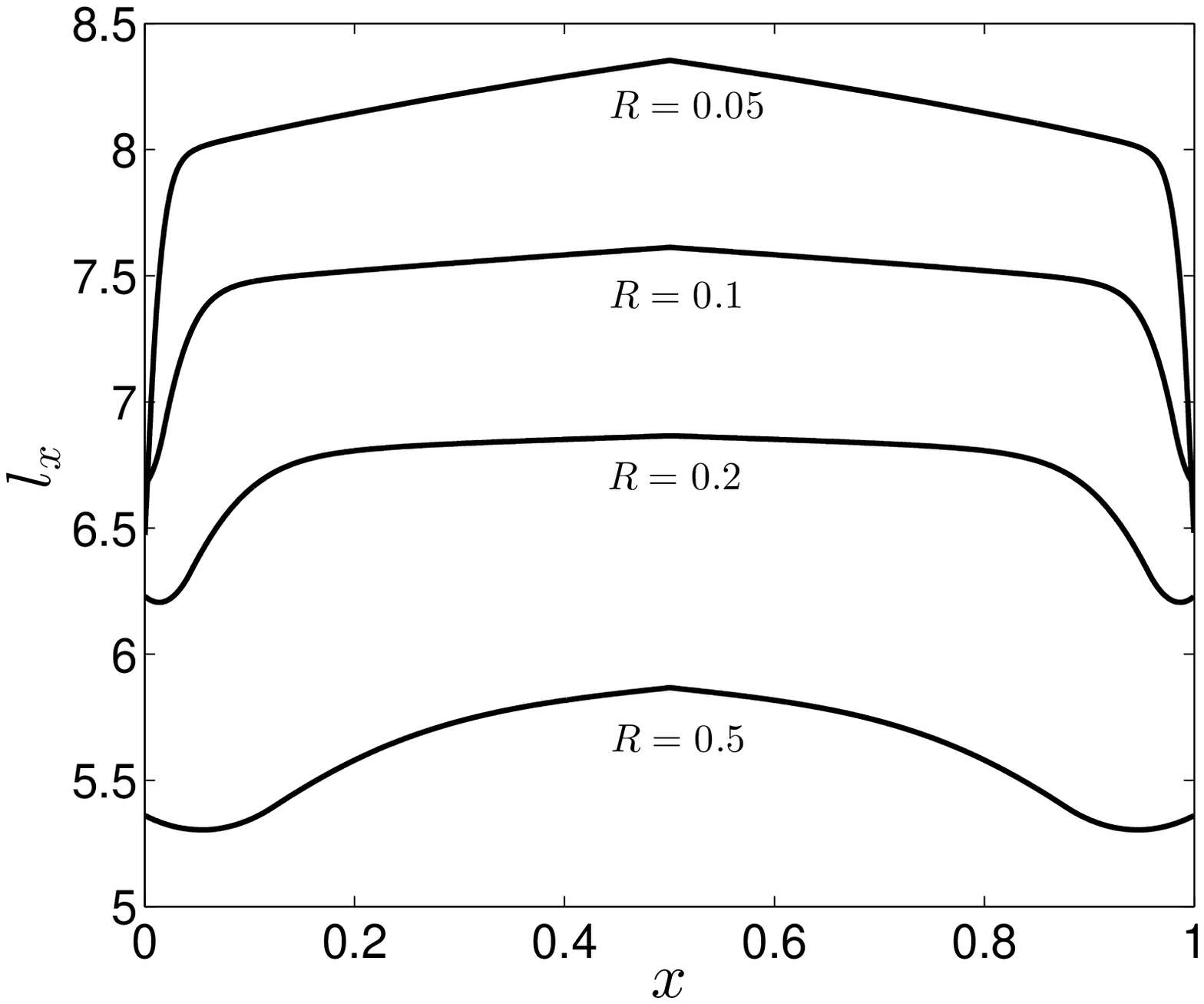}
\caption{$b=32$}
\label{fig:hfun_dc_1d-b32}
\end{subfigure}
\caption{ The excursion level functions $l_x$ corresponding to $b=4$ and $b=32$ for different correlation lengths $\cl$.
$f\equiv 1$.
}
\label{fig:hfun_dc_1d}
\end{center}
\end{figure}

\begin{figure}
\begin{center}
\begin{subfigure}[b]{0.48\textwidth}
\includegraphics[width=\textwidth]{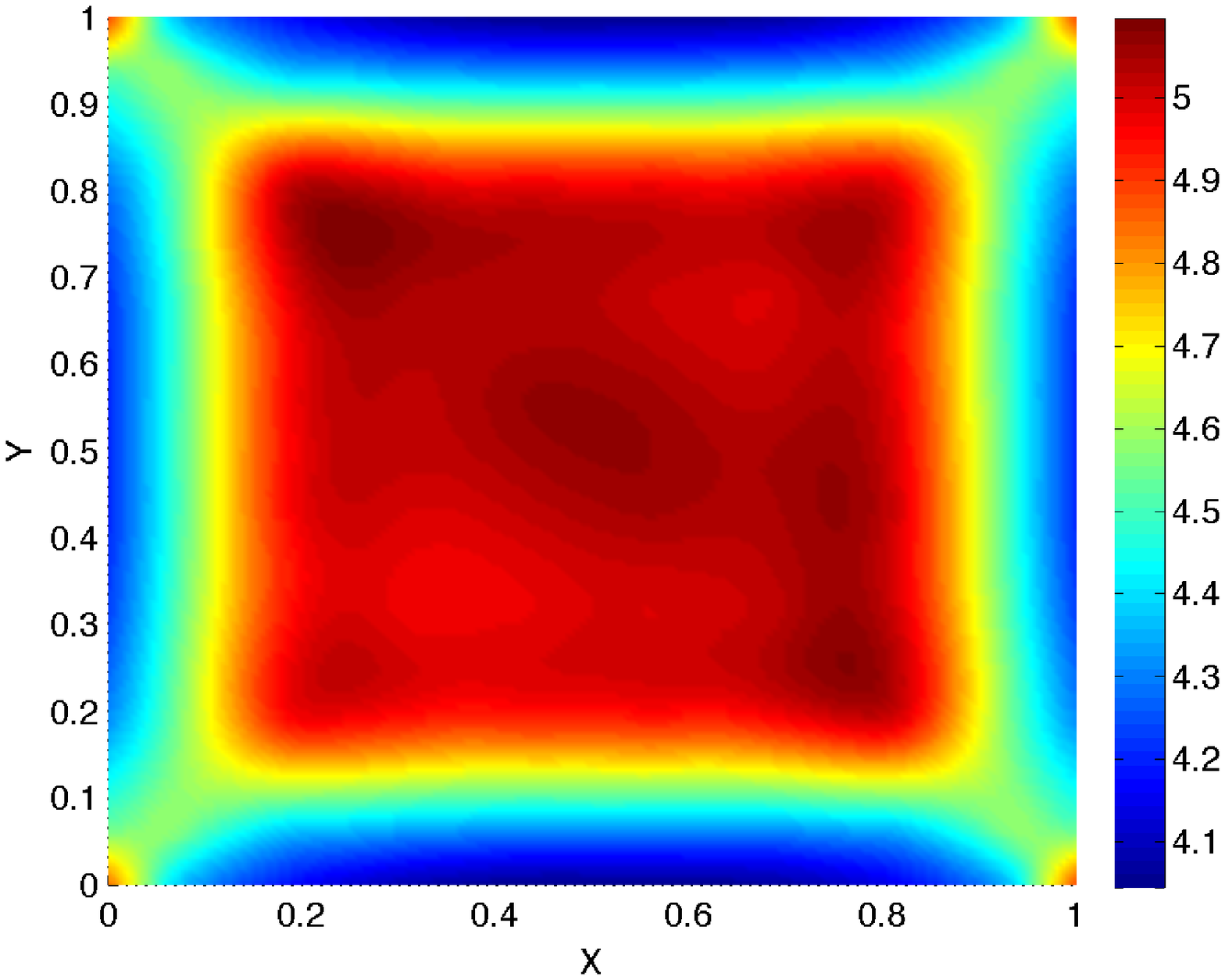}
\caption{$f\equiv 1$}
\label{figA:PeakHeight_b4_lc0_2_dU2_f1}
\end{subfigure}
\begin{subfigure}[b]{0.48\textwidth}
\includegraphics[width=\textwidth]{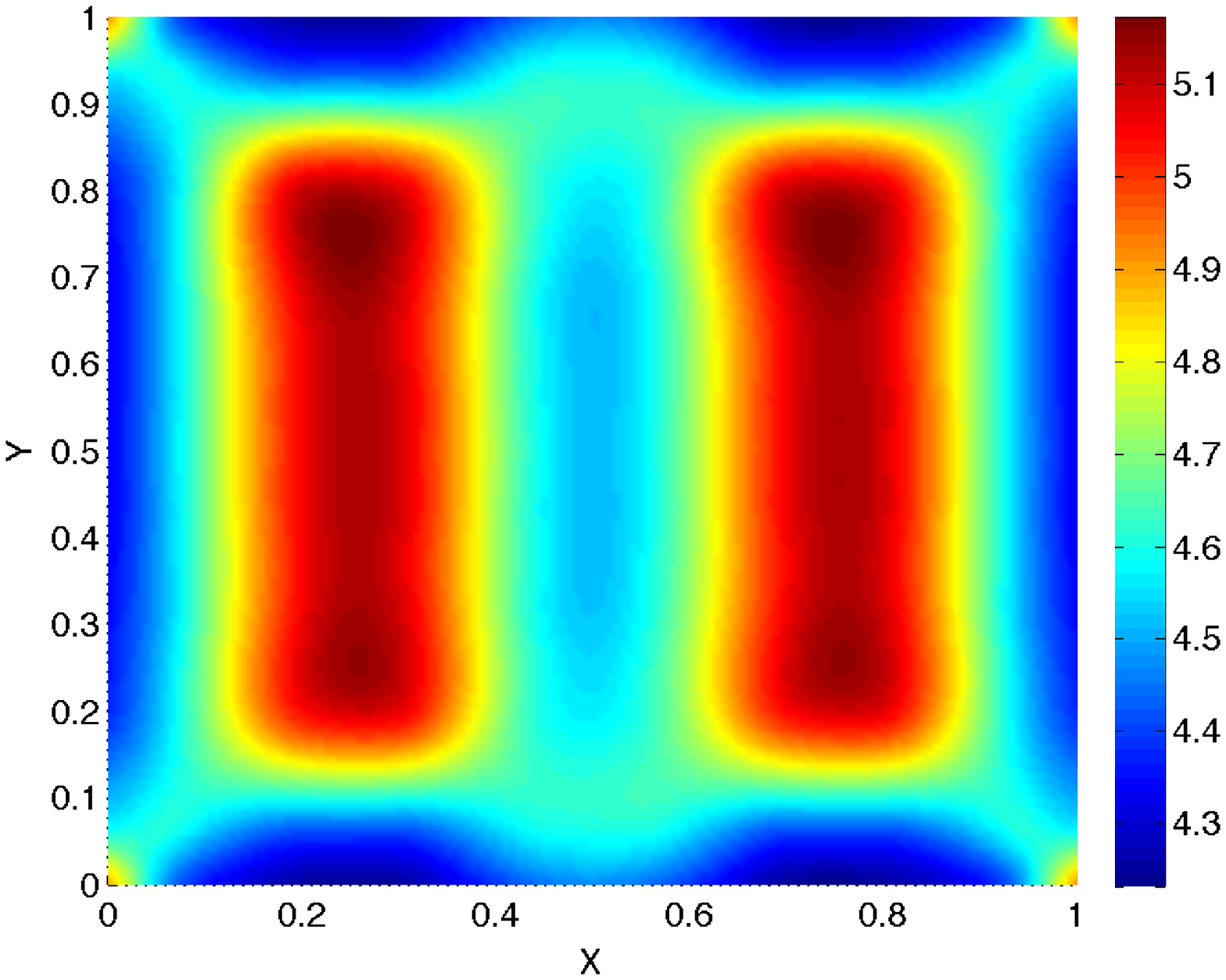}
\caption{{$f=\mbox{sign}( x_1\leq 0.5)$}}
\label{figB:PeakHeight_b4_lc0_2_dU2_f1}
\end{subfigure}
\caption{Contour plot of the excursion level function $l_x$ for reaching $\max |\nabla u|>b$ with homogeneous external force $f\equiv 1$ (left panel)
and the discontinuous external force $f=\mbox{sign}( x<0.5)$ (right panel).
The correlation length is  $\cl = 0.2$. The threshold is $b=4$.   }
\label{fig:PeakHeight_b4_lc0_2_dU2_f1}
\end{center}
\end{figure}

We now explore the function $l_x$ for different external force
functions $f$.  \figref{fig:hfun_dc_1d} shows the excursion level
functions for the one-dimensional problem with $\dom = [0,1]$ and $f$
is set to be a constant.  Due to the Dirichlet boundary condition, the
excursion level function $l_x$ is not a constant and it has
significantly lower values close to the boundary, especially when the
correlation length $\cl$ is small.  Thus, $\xi$ has a higher
probability to exhibit a high excursion near the boundary; so it is
for $\nabla u$.  The calculation of $l_x$ for the two-dimensional case
also confirms this boundary effect. Refer to \figref
{fig:PeakHeight_b4_lc0_2_dU2_f1} for the two choices of external
force.  In \figref {figA:PeakHeight_b4_lc0_2_dU2_f1}, the external
force is homogeneous $f(x) \equiv 1$.  In \figref
{figB:PeakHeight_b4_lc0_2_dU2_f1}, the external force has
discontinuity at $x_1=0.5$:
$$f(x) = \begin{cases}  1 &\mbox{if } x_1 \leq 0.5 \\
-1 & \mbox{if }  x_1>0.5. \end{cases} .$$
The local dip of the excursion level function $l_x$ near this discontinuity is consistent with the physical heuristics
that the material is easy to break down at this discontinuity line.

The dependence of the excursion level function on the correlation length $\cl$ is suggested in \figref{fig:hfun_dc_1d} from which smaller $\cl$ requires  larger $l_x$.
However, it should be noted that this does not imply that a smaller failure probability for smaller $\cl$  because it is easier for a Gaussian random function  with smaller $\cl$ to generate high excursions.

We check the effectiveness of the obtained excursion level function
for the importance sampling scheme by investigating the conditional
sample of $\xi$ given that the failure event $\sup|\nabla u|>b$
occurs.  Using direct Monte Carlo for a moderate $b$ with $f$ being
constant, a few samples were generated given the failure events.  We
observed a common feature of their spatial profiles from these
samples, each of which admits a unique and very high global maximum
for each sample of $\xi$ near the boundary, where the excursion level
function $l_x$ has dips.  It is worth mentioning that there are
several {\it local} maxima in the domain, but these local maxima are
significantly lower than the global one. Most of the strain fields
$|\nabla u|$ corresponding to these realizations of $\xi$ show a
global maximum close to the global maximizer of $\xi$.

\subsection{On the computation of the failure probabilities}

Table \ref{tab:lc01du1Da} summarises the performance of direct Monte
Carlo method (MC) and the proposed importance sampling method (IS)
described in Algorithm \ref{algdis} for the failure probabilities for
the one-dimensional differential equation.  Table \ref{tab:lc06du}
shows the results for two-dimensional case with two resolutions of
grid mesh, $25\times 25$ and $50\times 50$.  In these tables, ``$\wh
p_b$" columns are the estimated probabilities. The column ``std"
includes the standard deviation of one sample.  The ``rel. err."  is
the ratio of ``std" over ``$\wh p_b$".  For very large values of $b$,
{\it e.g.}, $b= 16$ or $32$, due to the finite number of samples,
direct Monte Carlo fails to yield reasonable estimates as the failure
event has not been observed in the samples. In this case, $\wh p_b$
and ``std" are marked as ``$-$" and the relative error, ``rel. err.",
is calculated by the theoretic result $\sqrt{1/p-1}$ with $p$ being
the estimated $\wh p_b$ from the importance sampling method.

%

\begin{table}[htdp]
\begin{center}
\begin{tabular}{|r|cc|cc|cc|}
\hline
&\multicolumn{2}{c|}{$\wh p_b$} &\multicolumn{2}{c|}{std}
&\multicolumn{2}{c|}{rel. err.}
\\
\hline
$b$  & MC
 &IS& MC & IS  &   MC & IS   \\
2 &2.15e-1& 2.15e-1
& 4.12e-1& 5.33-1
& 1.91& 2.48 \\
4 & 3.06e-2&3.02e-2&
1.74e-1&7.78e-2&
5.63&2.54\\
8 &2.75e-4&3.47e-4
&1.66e-2&8.29e-4
&60.4&2.39\\
16 &-&1.08e-5
&-&2.90e-5
&302*&2.69
\\
32 &-&1.89e-7&-&5.61e-7
&2300*&2.97 \\
\hline
\end{tabular}
\end{center}
*:  {\it\footnotesize relative error is given by $\sqrt{p^{-1}-1}$
  where $p$ is estimated from our proposed method.}

\caption{ (1D) The estimated failure probabilities $P(\max |\nabla u|>b)$ for the one-dimensional equation where $f\equiv 1$, $\cl=0.1$, and $N=400$ based on $10^6$ independent Monte Carlo samples in both direct Monte Carlo and the importance sampling.
}
\label{tab:lc01du1Da}
\end{table}%

\begin{table}[htdp]
\centering
\begin{subtable}[b]{0.8\textwidth}
\begin{center}
\begin{tabular}{|r|cc|cc|ll|}
\hline
&\multicolumn{2}{c|}{$\wh p_b$} &\multicolumn{2}{c|}{std}
&\multicolumn{2}{c|}{rel. err.}
\\
\hline
$b$  & MC
 &IS& MC & IS  &   MC & IS   \\
1&2.70e-1& 2.69e-1
& 4.44e-1& 2.53e-1
&1.65 & 0.94\\
2 &5.79e-2& 5.78e-2
& 2.34e-1& 6.16e-2
&4.03& 1.07\\
4 &6.25e-3&6.27e-3&
7.88e-2&7.53e-3&
12.6&1.20\\
8 &3.52e-4&3.57e-4
&1.88e-2&5.07e-4
&53.3&1.42\\
16 &8.00e-6&1.11e-5
&2.83e-3&1.82e-5
&353.0&1.64
\\
32 &-&1.96e-7&-&3.60e-7
&2253*&1.84 \\
\hline
\end{tabular}
\end{center}
\caption{mesh size $25\times25$}
\end{subtable}

\begin{subtable}[b]{0.8\textwidth}
\centering
\begin{center}
\begin{tabular}{|r|cc|cc|ll|}
\hline
&\multicolumn{2}{c|}{$\wh p_b$} &\multicolumn{2}{c|}{std}
&\multicolumn{2}{c|}{rel. err.}
\\
\hline
$b$  & MC
 &IS& MC & IS  &   MC & IS   \\
1&2.97e-1& 2.97e-1
& 4.57e-1& 2.80e-1
&1.54 & 0.94\\
2 &6.71e-2& 6.73e-2
& 2.50e-1& 7.26e-2
&3.73& 1.08\\
4 &7.75e-3&7.72e-3&
8.77e-2&9.27e-3&
11.3&1.20\\
8 &4.48e-4&4.66e-4
&2.12e-2&6.73e-4
&47.2&1.44\\
16 &1.80e-5&1.55e-5
&4.24e-3&2.56e-5
&236.0&1.64
\\
32 &-&2.93e-7&-&5.20e-7
&1847*&1.78\\
\hline
\end{tabular}
\end{center}
\caption{mesh size $50\times 50$}
\end{subtable}
\caption{ (2D) The estimated failure probabilities  $P(\max |\nabla u|>b)$ for the two-dimensional equation, where  $\cl=0.6$ and $f\equiv 1$.
 Sample size is $10^6$.  
 }\label{tab:lc06du}
\end{table}%

The relative error measures the relative efficiency of the Monte Carlo
schemes.  The comparison of relative errors in the last columns of
Table \ref{tab:lc01du1Da} and Table \ref{tab:lc06du} shows that, for
all values of the threshold $b$, the proposed importance sampling
scheme substantially outperform direct Monte Carlo.  When the event
becomes rarer, its advantage becomes more significant.  The importance
sampling scheme maintains a very mild increment of the relative error
that remains to be single digit even when the probability is as small
as $10^{-7}$.  Results based on different mesh sizes in Table
\ref{tab:lc06du} shows a relative difference around $10\%$. This
indicates that the spatial resolution is fine enough to get a reasonable accurate numerically obtained efficiency of the estimators.

In Algorithm \ref{algdis}, the alternative distribution
$g_\tau(\cdot)$ of the random variable $\xi(\tau)$ is suggested to the
Gaussian $\mathcal{N}(l_\tau, \sigma_\tau^2)$ and the variance of this
Gaussian is set to be $\sigma_\tau \sim ({l_\tau})^{-1}$.  To justify
this choice of the variance, we compare different constant values of
$\sigma_x\equiv \sigma$ and present the effect of $\sigma$ on the
performance of the resulting importance sampling scheme.  The
comparison is presented in Table \ref{tab:lc06sigma}.  The
(one-sample) standard deviation of the resulting importance sampling
scheme are for different $\sigma$ values. In addition, the typical
values of the reciprocal of the excursion level $l_x$ ($ x\in \dom$)
is also calculated for comparison.  As can be clearly seen from the
results in the table, the optimal choice of $\sigma_\tau$ is indeed
around the reciprocal of the excursion level, $\frac{1}{l_\tau}$.

\begin{table}[htbp]
\begin{center}
\begin{tabular}{|c|c||c|c|c|c|c|c|c|}
\hline
$b$& $1/l_x$ &$\sigma=0.02$ & $\sigma=0.1$& $\sigma=0.2$ & $\sigma=0.3$& $\sigma=0.5$ & $\sigma=1$ & $\sigma=5$
\\
\hline
4  & $0.24\sim0.30$ &7.40e-1  & 8.67e-3& 9.43e-3   & 7.36e-3&8.66e-3&1.21e-2&2.81e-2\\
32  & $0.16\sim 0.18$ &1.65e-6 & 3.07e-7& 2.85e-7  &3.22e-7&4.05e-7&5.61e-7&1.40e-6\\
\hline
\end{tabular}
\end{center}
\caption{The standard deviation  of the importance sampling estimator for difference choices of the variance  $\sigma_\tau^2\equiv \sigma^2$ when conditionally sampling $\xi(\tau)$.}
\label{tab:lc06sigma}
\end{table}

For a smaller correlation length $\cl=0.2$, we test the algorithm with constant external force $f$. As the correlation length of the random field is smaller, we use a finer mesh grid ($150\times150$) to resolve. 
The results obtained are shown in Table \ref{tab:lc02du}, which further confirms the efficiency of the importance sampling scheme.

\begin{table}[htdp]
\begin{center}
\begin{tabular}{|r|cc|cc|ll|}
\hline
&\multicolumn{2}{c|}{$\wh p_b$} &\multicolumn{2}{c|}{std}
&\multicolumn{2}{c|}{rel. err.}
\\
\hline
$b$ & MC
 &IS& MC & IS  &   MC & IS   \\
1& 3.58e-1&  3.58e-1
& 4.80e-1  & 4.60e-1
& 1.34& 1.29 \\
2&3.02e-2 & 2.96e-2
&1.71e-2  & 1.67e-2
& 5.67 &5.64 \\
4& 8.65e-4 & 8.61e-4
& 2.94e-2 & 5.81e-3
& 34.0 &  6.75 \\
8&1.20e-5 &  1.44e-5
& 3.46e-3 & 6.28e-5
& 289& 4.36 \\
16&- & 1.54e-7
& - & 6.73e-7
& 2548*& 4.37 \\
32&- &  1.06e-9
& - & 1.29e-8
&30715* &  12.1\\
\hline
\end{tabular}
\end{center}
\caption
{(2D) $P(\max |\nabla u|>b)$.  Correlation length $\cl=0.2$.  Mesh size $150\times150$. Sample size $10^6$.
$f\equiv 1$.}
\label{tab:lc02du}
\end{table}%

\subsection{On the asymptotics of the failure probabilities}
\begin{figure}[tbph]
\begin{center}
\includegraphics[scale=0.46]{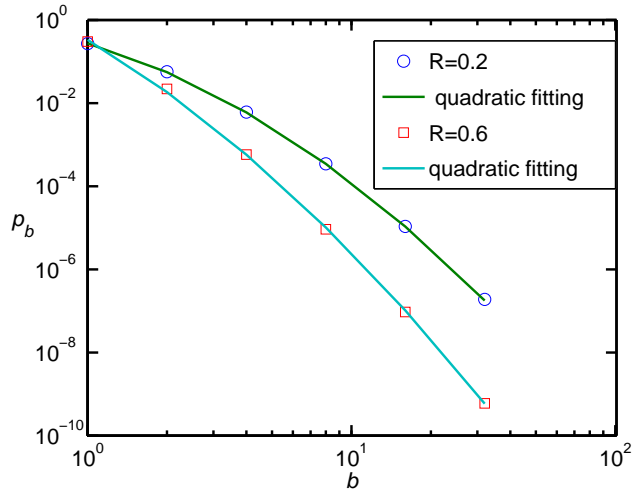}
\caption{log-log plot of $p_b$ vs $b$
from the data in Table \ref{tab:lc06du}(B) and  Table \ref{tab:lc02du}.
Set
$\log p_b = {q(\log b)}$ and use the quadratic function for $q$
for least square fitting.
}
\label{fig:2Db_pb_loglogfit}
\end{center}
\end{figure}

The importance sampling method can be applied to efficiently calculate
quantities related to the failure event.  A direct application is that
we can numerically characterize the asymptotic behavior of tail
probabilities $p_b=P(\sup_{x\in \dom}|\nabla u(x)| > b )$.  For
example, the data in Table \ref{tab:lc06du}(B) and Table
\ref{tab:lc02du} allow us to postulate an empirical asymptotics
between $p_b=P(\max|\nabla u|>b)$ and $b$.
\figref{fig:2Db_pb_loglogfit} shows the log-log plot of $p_b$ vs $b$
and the result of least square fitting. The result shows that the tail
distribution satisfies $$p_b = e^{ q(\log b) },$$ where $q$ is a
quadratic function and refer to \figref{fig:2Db_pb_loglogfit} for
specific expression and its dependency on the correlation length
$\cl$.  This quadratic dependency is consistent with our analytical
result for one dimensional case \cites{LiuZhou13,LiuZhou14}.

In addition, Table \ref{tab:lc06du} and Table \ref{tab:lc02du}
together suggest that smaller correlation length leads to larger
failure probabilities.  This observation is the same  with  the higher
excursion probabilities of Gaussian random fields with smaller
correlation length, although the random function $|\nabla u|$ is not 
trivially Gaussian.

\subsection{Numerical results under periodic boundary condition}
We have demonstrated the efficiency of the variance reduction on our
importance sampling scheme for the Dirichlet boundary condition and
the constant external force.  Here, we test our method for the
periodic boundary condition and non homogeneous body force.  We only
study a one-dimensional example where the domain of the elliptic
equation is $\dom=[0,1]$.  To accommodate the periodic boundary
condition, we first sample the values of the periodic random field
$\xi(x)$ over the finite domain $[0,1]$.  This is done by simply
designating a periodic covariance function $C_p(x)$.  The function
$C_p$ we choose is a period-$1$ extension of the original covariance
function $C(x)=e^{-x^2/\cl^2}$ where $C_p(x)=C(x)$ for all $x\in
[-1/2,1/2]$. The solvability of the equation requires that $\int_0^1
f(x)=0$.  We consider the following quadratic function
$f=10(x-1/2)^2-5/6$ in our example.  The performance of the importance
sampling scheme is shown in Table \ref{tab:lc01du1D_PBC}.  The
variance reduction is not as significant as the Dirichlet condition
case, but still much better than the direct Monte Carlo for large $b$.

%
%
%
%

\begin{table}[ht]
\begin{center}
\begin{tabular}{|r|cc|cc|ll|}
\hline
&\multicolumn{2}{c|}{$\wh p_b$} &\multicolumn{2}{c|}{std}
&\multicolumn{2}{c|}{rel. err.}
\\
\hline
$b$  & MC
 &IS& MC & IS  &   MC & IS   \\
1 & 3.77e-2&3.76e-2&
1.90e-1&2.28e-1&
5.05&6.06\\
2 & 2.76e-3& 2.40e-3
& 5.25-2&8.80e-2
& 19.0& 36.7\\
4 & 9.68e-5& 7.17e-5
& 9.84-3& 3.22-3
&102 & 45.9
\\
8 &1.00e-6& 1.51e-6&
1.00e-3& 1.82e-4
& 1000& 121 \\
12 &- & 1.28e-7  &
-&1.77e-5
&2795* &138
\\
\hline
\end{tabular}

\end{center}
\caption{ (1D) $P(\max |\nabla u|>b)$.  Periodic boundary condition.  $\cl=0.2$. $Nx=400$.
Sample Size $4\times 10^6$.
$f=  10(x-0.5)^2-5/6$.
}\label{tab:lc01du1D_PBC}

\end{table}%

\section{Conclusion and discussion}\label{sec:conclude}

We present in this work an efficient importance sampling strategy for
computing small probabilities associated with materials failure
problem modeled by a scalar elliptic equation with random log-normal
coefficient. The change of measure used in the importance sampling is
suggested by one dimensional analysis and further justified
numerically for higher dimensions. Our numerical results verifies the
superior behavior of the estimator over conventional approaches.

The asymptotic analysis of the failure probability, in particular,
rigorously establishing the relation \eqref{eqn:exp-scale} is a very
interesting open problem.

In this work, for simplicity, we have used the scalar model, however,
there is no conceptual difficulty in generalizing our method to linear
elastic models, which would be of interest for practical applications.

The efficient sampling technique for the failure event opens the door
to many interesting applications. One particularly interesting
application would be the design of materials in consideration of
minimizing the failure probability. This gives arise to interesting
future directions to explore.

\section*{Acknowledgments}
J. Liu is supported in part by NSF CMMI-1069064, NSF SES-1323977, and
Army Grant W911NF-14-1-0020.  The work of J. Lu was supported in part
by the Alfred P.~Sloan Foundation and the NSF grant DMS-1312659.
X. Zhou acknowledges the financial support from CityU Start-Up Grant
(7200301) and Hong Kong Early Career Schemes (109113).

\bibliographystyle{amsxport}
\bibliography{MaterialsFailure}

\end{document}